\documentclass[12pt,reqno]{amsart}
\usepackage{amsmath,amssymb,amsthm,mathtools,enumerate}
\DeclareFontEncoding{OT2}{}{}
\DeclareFontSubstitution{OT2}{cmr}{m}{n}
\DeclareFontFamily{OT2}{cmr}{\hyphenchar\font45 }
\DeclareFontShape{OT2}{cmr}{m}{n}{
<5><6><7><8><9>gen*wncyr
<10><10.95><12><14.4><17.28><20.74><24.88>wncyr10}{}
\DeclareFontShape{OT2}{cmr}{b}{n}{
<5><6><7><8><9>gen*wncyb
<10><10.95><12><14.4><17.28><20.74><24.88>wncyb10}{}
\DeclareMathAlphabet{\mathcyr}{OT2}{cmr}{m}{n}
\DeclareMathAlphabet{\mathcyb}{OT2}{cmr}{b}{n}
\SetMathAlphabet{\mathcyr}{bold}{OT2}{cmr}{b}{n}

\calclayout
\theoremstyle{plain}
\newtheorem{thm}{Theorem}[section]

\newtheorem{proposition}[thm]{Proposition}
\theoremstyle{definition}

\theoremstyle{remark}

\newcommand{\bk}{\boldsymbol{k}} 
\newcommand{\bl}{\boldsymbol{l}}
\newcommand{\bh}{\boldsymbol{h}}
\newcommand{\bn}{\boldsymbol{n}}
\newcommand{\bm}{\boldsymbol{m}}
\newcommand{\br}{\boldsymbol{r}}
\newcommand{\wt}{\mathrm{wt}}
\newcommand{\ZZ}{\mathbb{Z}}
\newcommand{\II}{\mathcal{I}}
\newcommand{\sh}{\mathbin{\mathcyr{sh}}}
\title{Connectors}
\author{Shin-ichiro Seki}
\address{Mathematical Institute, Tohoku University, 6-3, Aoba, Aramaki, Aoba-Ku,
Sendai, 980-8578, Japan}
\email{shinichiro.seki.b3@tohoku.ac.jp}
\begin{document}
\begin{abstract}
Recently, the author and Yamamoto invented a new proof of the duality for multiple zeta values.
The technique is applicable in other series identities.
In this article, we exhibit such proofs for some series identities.
\end{abstract}
\maketitle
\setcounter{tocdepth}{1}
\tableofcontents
\section{Dynamic (or algorithmic) proofs}
For a given series identity, the strategy to prove it in this article is as follows:

\

\begin{enumerate}[(i)]
	\item Define the \emph{connected sum} with the \emph{connector} or the \emph{connecting relation}.
	\item\label{sec1:2} Check its \emph{transport relations}.
	\item Transport indices, algorithmically.
	\item\label{sec1:4} Check the \emph{boundary conditions}.
\end{enumerate}

\

Only these!

\

Proofs of \eqref{sec1:2} and \eqref{sec1:4} are relatively easier than the proof of the original identity itself.
We can often prove them by using partial fraction decompositions or telescoping sums and we leave the proofs to the reader.
The hardest part in each dynamic proof is how to find a suitable connected sum.

Note that values of connected sums may diverge, however we use only convergent values when we apply them.
\section{Notation}
$N$ denotes a positive integer.
$p$ denotes a prime number.
For a positive integer $a$, we set
\[
J_a\coloneqq\underbrace{\ZZ_{\geq0}\times\cdots\times\ZZ_{\geq0}}_a,\quad I_a\coloneqq\underbrace{\ZZ_{\geq1}\times\cdots\times\ZZ_{\geq1}}_a,
\]
\[
I'_a\coloneqq\{(k_1,\dots,k_a)\in I_a : k_a\geq 2\}.
\]
Further, $J_0=I_0=I'_0\coloneqq\{\varnothing\}$ and
\[
\II\coloneqq\bigsqcup_{a=1}^{\infty}I_a,\quad \II'\coloneqq\bigsqcup_{a=1}^{\infty}I'_a,\quad\II_0\coloneqq\bigsqcup_{a=0}^{\infty}I_a,\quad \II_0'\coloneqq\bigsqcup_{a=0}^{\infty}I'_a.
\]

Let $a$ be a positive integer and $\bk=(k_1,\dots,k_a)\in I_a$.
The weight of $\bk$ is defined as $\wt(\bk)\coloneqq k_1+\cdots+k_a$. 
The reverse index of $\bk$ is defined as $\overline{\bk}=(k_a,\dots,k_1)$.
For $0\leq i\leq a$, $\bk_{(i)}\coloneqq(k_1,\dots,k_i)$, $\bk^{(i)}\coloneqq(k_{i+1},\dots,k_a)$ ($\bk_{(0)}=\bk^{(a)}=\varnothing$).
Further, we define the \emph{arrow-notation} as
\begin{align*}
\bk_{\to}&\coloneqq(k_1,\dots,k_a,1),\qquad\qquad{}_{\leftarrow}\bk\coloneqq(1,k_1,\dots,k_a),\\
\bk_{\uparrow}&\coloneqq(k_1,\dots,k_{a-1},k_a+1),\quad{}_{\uparrow}\bk\coloneqq(k_1+1,k_2,\dots,k_a),\\
\bk_{\downarrow}&\coloneqq(k_1,\dots,k_{a-1},k_a-1),\quad{}_{\downarrow}\bk\coloneqq(k_1-1,k_2,\dots,k_a).
\end{align*}
$\{\uparrow\}^j$ denotes $\underbrace{\uparrow\cdots\uparrow}_j$ and $\{\downarrow\}^j$ denotes $\underbrace{\downarrow\cdots\downarrow}_j$.
By convention, we use $\wt(\varnothing)=0$, $\overline{\varnothing}=\varnothing$, $\varnothing_{\to}={}_{\leftarrow}\varnothing=(1)$, and $\varnothing_{\uparrow}={}_{\uparrow}\varnothing=\varnothing_{\downarrow}={}_{\downarrow}\varnothing=\varnothing$.

For $\bk=(k_1,\dots,k_a)\in I_a$, $\bl=(l_1,\dots,l_b)\in I_b$ ($a,b\geq 0$), $(\bk,\bl)$ denotes the concatenation index $(k_1,\dots, k_a, l_1,\dots, l_b)$.

In definitions of connected sums and \eqref{eq:MHS}, we use abbreviated notation $\bn,\bm,\br\in\II_0$ as
\[
\bn=(n_1,\dots,n_a),\quad \bm=(m_1,\dots,m_b),\quad \br=(r_1,\dots,r_c),
\]
even if there appear extra variables $n_0$, $m_{b+1}$ and so on.

For $\bk=(k_1,\dots,k_a)\in J_a$ and $\bn=(n_1,\dots, n_a)\in I_a$, $\bn^{\bk}\coloneqq\prod_{i=1}^an_i^{k_i}$ ($\bn^{\bk}\coloneqq1$ when $a=0$).

For $\bk\in I_a$, we define $\zeta^{}_N(\bk)$ and $\zeta^{\star}_N(\bk)$ as
\begin{equation}\label{eq:MHS}
\zeta^{}_N(\bk)\coloneqq\sum_{0=n_0< n_1<\cdots<n_a\leq N}\frac{1}{\bn^{\bk}},\quad \zeta^{\star}_N(\bk)\coloneqq\sum_{0=n_0< n_1\leq\cdots\leq n_a\leq N}\frac{1}{\bn^{\bk}}.
\end{equation}
For $\bk\in \II'_0$, $\zeta(\bk)\coloneqq\lim\limits_{N\to\infty}\zeta^{}_N(\bk)<+\infty$.
For a connected sum $Z_N$, $Z_{\infty}$ means $\lim\limits_{N\to\infty}Z_N$.

We omit to define some standard notions including the shuffle product $\sh$, the harmonic product $\ast$, the dual index $\bk^{\dagger}$, and the Hoffman dual index $\bk^{\vee}$.
Rather, we can define these notions by the following transport relations and algorithms.
\section{Shuffle product formula}
\begin{proposition}[{\cite[Theorem~4.1]{HO}}]
For $\bk,\bl\in \II'_0$, we have
\begin{equation}\label{eq:shuffle1}
\zeta(\bk)\zeta(\bl)=\sum a_{\bh}\zeta(\bh),
\end{equation}
where $\bh\in \II'_0$ runs over indices appearing in $\bk\sh\bl=\sum a_{\bh}\bh$.
\end{proposition}
\begin{proposition}[\cite{KZ}, {\cite[Corollary~4.1]{On}}]
For $\bk,\bl\in \II_0$, we have
\begin{equation}\label{eq:shuffle1'}
(-1)^{\wt(\bl)}\zeta^{}_{p-1}(\bk,\overline{\bl})\equiv\sum a_{\bh}\zeta^{}_{p-1}(\bh)\pmod{p},
\end{equation}
where $\bh\in \II_0$ runs over indices appearing in $\bk\sh\bl=\sum a_{\bh}\bh$.
\end{proposition}
\subsection{Connected sum}
For $\bk\in I_a, \bl\in I_b, \bh\in I_c$ ($a,b\in\ZZ_{\geq 0}, c\in \ZZ_{\geq 1}$), the connected sum $Z_N^{\sh}(\bk;\bl;\bh)$ is defined by
\[    
Z_N^{\sh}(\bk;\bl;\bh)\coloneqq\sum_{\substack{0=n_0<n_1<\cdots<n_a \\ 0=m_0<m_1<\cdots<m_b \\ n_a+m_b=r_1<r_2<\cdots<r_c\leq N}}\frac{1}{\bn^{(\bk_{\downarrow})}\bm^{(\bl_{\downarrow})}\br^{({}_{\downarrow}\bh)}}.
\]
The connecting relation is $n_a+m_b=r_1$.
By definition, there is a symmetry
\begin{equation}\label{eq:shuffle2}
Z_N^{\sh}(\bk;\bl;\bh)=Z_N^{\sh}(\bl;\bk;\bh).
\end{equation}
This connected sum is essentially defined in the right-hand side of \cite[(18)]{KMT}.
\subsection{Transport relations}
For $\bk,\bl\in\II_0, \bh\in \II$, we have
\begin{equation}\label{eq:shuffle3}
Z_N^{\sh}(\bk_{\uparrow};\bl_{\uparrow};\bh)=Z_N^{\sh}(\bk;\bl_{\uparrow};{}_{\uparrow}\bh)+Z_N^{\sh}(\bk_{\uparrow};\bl;{}_{\uparrow}\bh)\quad (\bk,\bl\neq\varnothing),
\end{equation}
\begin{equation}\label{eq:shuffle4}
Z_N^{\sh}(\bk_{\to};\bl;{}_{\uparrow}\bh)=Z_N^{\sh}(\bk_{\uparrow};\bl;{}_{\leftarrow}\bh).
\end{equation}
\subsection{Algorithm}\label{sh-alg}
For each $Z_N^{\sh}(\bk;\bl;\bh)$, we set rules as follows:

\

\begin{enumerate}[(i)]
\item\label{sec3:1} If $\bk$ and $\bl\in\II'$, then use the transport relation \eqref{eq:shuffle3}.
\item\label{sec3:2} If $\bl\in\II_0\setminus\II'$, then use the symmetry \eqref{eq:shuffle2}.
\item\label{sec3:3} If $\bk\in\II\setminus\II'$, then use the transport relation \eqref{eq:shuffle4}.
\item\label{sec3:4} If $\bk=\varnothing$, then stop.
\end{enumerate}

\

Start from $Z_N^{\sh}(\bk_{\uparrow};\bl_{\uparrow};(1))$ ($\bk,\bl\in\II_0$) and transport indices according to the above rules until the algorithm stops for all connected sums. Then, we have an identity which has the following form ($\bh,\bh'\in\II_0$):
\begin{equation}\label{eq:shuffle5}
Z_N^{\sh}(\bk_{\uparrow};\bl_{\uparrow};(1))=\sum Z_N^{\sh}(\varnothing;\bh_{\uparrow};{}_{\leftarrow}\bh').
\end{equation}
\subsection{Boundary conditions}
For $\bk,\bl\in\II'_0$, we have
\begin{equation}\label{eq:shuffle6}
Z_{\infty}^{\sh}(\bk_{\uparrow};\bl_{\uparrow};(1))=\zeta(\bk)\zeta(\bl)
\end{equation}
and 
\[
Z_N^{\sh}(\varnothing;\bh_{\uparrow};{}_{\leftarrow}\bh')=\zeta^{}_N(\bh,\bh').
\]
By these boundary conditions, we see that the case $N\to\infty$ of the identity \eqref{eq:shuffle5} is nothing but the shuffle product formula \eqref{eq:shuffle1}.

When $\bk,\bl\in\II_0$, by using the following congruence instead of using the boundary condition \eqref{eq:shuffle6}, we can also prove the shuffle relation for finite multiple zeta values \eqref{eq:shuffle1'}:
\[
Z_{p-1}^{\sh}(\bk_{\uparrow};\bl_{\uparrow};(1))\equiv (-1)^{\wt(\bl)}\zeta^{}_{p-1}(\bk,\overline{\bl})\pmod{p}.
\]
\subsection{Example}
By applying the algorithm in \ref{sh-alg} to $Z_N^{\sh}((1,2);(2);(1))$, we have
\begin{align*}
&Z_N^{\sh}((1,2);(2);(1))\stackrel{\eqref{sec3:1}}{=}Z_N^{\sh}((1,1);(2);(2))+Z_N^{\sh}((1,2);(1);(2)),\\
&\qquad Z_N^{\sh}((1,1);(2);(2))\stackrel{\eqref{sec3:3}}{=}Z_N^{\sh}((2);(2);(1,1))\stackrel{\eqref{sec3:1}}{=}Z_N^{\sh}((1);(2);(2,1))+Z_N^{\sh}((2);(1);(2,1)),\\
&\qquad\qquad Z_N^{\sh}((2);(1);(2,1))\stackrel{\eqref{sec3:2}}{=}Z_N^{\sh}((1);(2);(2,1))\stackrel{\eqref{sec3:3}}{=}Z_N^{\sh}(\varnothing;(2);(1,1,1)),\\
&\qquad Z_N^{\sh}((1,2);(1);(2))\stackrel{\eqref{sec3:2}}{=}Z_N^{\sh}((1);(1,2);(2))\stackrel{\eqref{sec3:3}}{=}Z_N^{\sh}(\varnothing;(1,2);(1,1)).
\end{align*}
Thus,
\[
Z_N^{\sh}((1,1)_{\uparrow};(1)_{\uparrow};(1))=2Z_N^{\sh}(\varnothing;(1)_{\uparrow};{}_{\leftarrow}(1,1))+Z_N^{\sh}(\varnothing;(1,1)_{\uparrow};{}_{\leftarrow}(1))
\]
and this corresponds to
\[
(1,1)\sh(1)=2((1),(1,1))+((1,1),(1))=3(1,1,1).
\]
\section{Harmonic product formula}
\begin{proposition}[{\cite[Theorem~3.2]{H2}}]
For $\bk,\bl\in \II_0$, we have
\begin{equation}\label{eq:harmonic1}
\zeta^{}_N(\bk)\zeta^{}_N(\bl)=\sum b_{\bh}\zeta^{}_N(\bh),
\end{equation}
where $\bh\in \II_0$ runs over indices appearing in $\bk\ast\bl=\sum b_{\bh}\bh$.
\end{proposition}
\subsection{Connected sum}
For $\bk\in I_a, \bl\in I_b, \bh\in I_c$ ($a,b\in\ZZ_{\geq 1}, c\in\ZZ_{\geq 0}$), the connected sum $Z_N^{\ast}(\bk;\bl;\bh)$ is defined by
\[    
Z_N^{\ast}(\bk;\bl;\bh)\coloneqq\sum_{\substack{1=r_0\leq r_1<\cdots<r_c \\ r_c=n_1<\cdots<n_a\leq N \\ r_c=m_1<\cdots<m_b\leq N}}\frac{1}{\bn^{({}_{\downarrow}\bk)}\bm^{({}_{\downarrow}\bl)}\br^{\bh}}.
\]
The connecting relation is $n_1=m_1=r_c$.
By definition, there is a symmetry
\begin{equation}\label{eq:harmonic2}
Z_N^{\ast}(\bk;\bl;\bh)=Z_N^{\ast}(\bl;\bk;\bh).
\end{equation}
The definition of this connected sum is essentially due to Hirose.
\subsection{Transport relations}
For $\bk,\bl\in\II, \bh\in \II_0$, we have
\begin{equation}\label{eq:harmonic3}
Z_N^{\ast}({}_{\leftarrow}\bk;{}_{\leftarrow}\bl;\bh)=
Z_N^{\ast}(\bk;{}_{\leftarrow}\bl;\bh_{\to})
+Z_N^{\ast}({}_{\leftarrow}\bk;\bl;\bh_{\to})
+Z_N^{\ast}(\bk;\bl;\bh_{\to\uparrow}),
\end{equation}
\begin{equation}\label{eq:harmonic4}
Z_N^{\ast}({}_{\uparrow}\bk;\bl;\bh)=Z_N^{\ast}(\bk;\bl;\bh_{\uparrow}).
\end{equation}
\subsection{Algorithm}\label{har-alg}
For each $Z_N^{\ast}(\bk;\bl;\bh)$, we set rules as follows:

\

\begin{enumerate}[(i)]
\item\label{sec4:1} If $\overline{\bk}$ and $\overline{\bl}\in\II\setminus(\II'\cup\{(1)\})$, then use the transport relation \eqref{eq:harmonic3}.
\item\label{sec4:2} If $\overline{\bl}\in\II'\cup\{( 1)\}$, then use the symmetry \eqref{eq:harmonic2}.
\item\label{sec4:3} If $\overline{\bk}\in\II'$, then use the transport relation \eqref{eq:harmonic4}.
\item\label{sec4:4} If $\bk=(1)$, then stop.
\end{enumerate}

\

Start from $Z_N^{\ast}({}_{\leftarrow}\bk;{}_{\leftarrow}\bl;\varnothing)$ ($\bk,\bl\in\II_0$) and transport indices according to the above rules until the algorithm stops for all connected sums. Then, we have an identity which has the following form ($\bh,\bh'\in\II_0, h\in\ZZ_{\geq 0}$):
\begin{equation}\label{eq:harmonic5}
Z_N^{\ast}({}_{\leftarrow}\bk;{}_{\leftarrow}\bl;\varnothing)=\sum Z_N^{\ast}((1);{}_{\{\uparrow\}^h\leftarrow}\bh;\bh').
\end{equation}
\subsection{Boundary conditions}
For $\bk,\bl,\bh,\bh'\in\II_0, h\in\ZZ_{\geq 0}$, we have
\[
Z_N^{\ast}({}_{\leftarrow}\bk;{}_{\leftarrow}\bl;\varnothing)=\zeta^{}_N(\bk)\zeta^{}_N(\bl)
\]
and
\[
Z_N^{\ast}((1);{}_{\{\uparrow\}^h\leftarrow}\bh;\bh')=\zeta^{}_N(\bh'_{\{\uparrow\}^{h}},\bh).
\]
By these boundary conditions, we see that the identity \eqref{eq:harmonic5} is nothing but the harmonic product formula \eqref{eq:harmonic1}.
\subsection{Example}
By applying the algorithm in \ref{har-alg} to $Z_N^{\ast}((1,1);(1,2);\varnothing)$, we have
\begin{align*}
&Z_N^{\ast}((1,1);(1,2);\varnothing)\stackrel{\eqref{sec4:1}}{=}Z_N^{\ast}((1);(1,2);(1))+Z_N^{\ast}((1,1);(2);(1))+Z_N^{\ast}((1);(2);(2)),\\
&\qquad Z_N^{\ast}((1,1);(2);(1))\stackrel{\eqref{sec4:2}}{=}Z_N^{\ast}((2);(1,1);(1))\stackrel{\eqref{sec4:3}}{=}Z_N^{\ast}((1);(1,1);(2)).
\end{align*}
Thus,
\[
Z_N^{\ast}({}_{\leftarrow}(1);{}_{\leftarrow}(2);\varnothing)=Z_N^{\ast}((1);{}_{\leftarrow}(2);(1))+Z_N^{\ast}((1);{}_{\leftarrow}(1);(2))+Z_N^{\ast}((1);{}_{\uparrow\leftarrow}\varnothing;(2))
\]
and this corresponds to
\[
(1)\ast(2)=(1,2)+(2,1)+(3).
\]

\

By applying the algorithm in \ref{har-alg} to $Z_N^{\ast}((1,1,1);(1,1);\varnothing)$, we have
\begin{align*}
&Z_N^{\ast}((1,1,1);(1,1);\varnothing)\stackrel{\eqref{sec4:1}}{=}Z_N^{\ast}((1,1);(1,1);(1))+Z_N^{\ast}((1,1,1);(1);(1))+Z_N^{\ast}((1,1);(1);(2)),\\
&\quad Z_N^{\ast}((1,1);(1,1);(1))\stackrel{\eqref{sec4:1}}{=}Z_N^{\ast}((1);(1,1);(1,1))+Z_N^{\ast}((1,1);(1);(1,1))+Z_N^{\ast}((1);(1);(1,2)),\\
&\quad\quad Z_N^{\ast}((1,1);(1);(1,1))\stackrel{\eqref{sec4:2}}{=}Z_N^{\ast}((1);(1,1);(1,1)),\\
&\quad Z_N^{\ast}((1,1,1);(1);(1))\stackrel{\eqref{sec4:2}}{=}Z_N^{\ast}((1);(1,1,1);(1)),\\
&\quad Z_N^{\ast}((1,1);(1);(2))\stackrel{\eqref{sec4:2}}{=}Z_N^{\ast}((1);(1,1);(2)).
\end{align*}
Thus,
\begin{align*}
Z_N^{\ast}({}_{\leftarrow}(1,1);{}_{\leftarrow}(1);\varnothing)&=2Z_N^{\ast}((1);{}_{\leftarrow}(1);(1,1))+Z_N^{\ast}((1);{}_{\leftarrow}\varnothing;(1,2))+Z_N^{\ast}((1);{}_{\leftarrow}(1,1);(1))\\
&\qquad+Z_N^{\ast}((1);{}_{\leftarrow}(1);(2))
\end{align*}
and this corresponds to
\[
(1,1)\ast(1)=2((1,1),(1))+(1,2)+((1),(1,1))+(2,1)=3(1,1,1)+(1,2)+(2,1).
\]
\section{Duality}
Let $n,m$ be non-negative integers.
For $\bk\in I'_a$, we define the \emph{multiple zeta value with double tails} $\zeta^{}_{n,m}(\bk)$ by
\[
\zeta^{}_{n,m}(\bk)\coloneqq\sum_{n=n_0< n_1<\cdots<n_a}\frac{1}{\bn^{\bk}}\cdot\frac{1}{\binom{n_a+m}{m}}.
\]
\begin{proposition}[\cite{A}]
For any $\bk\in\II'$,
\begin{equation}\label{eq:duality1}
\zeta^{}_{n,m}(\bk)=\zeta^{}_{m,n}(\bk^{\dagger}).
\end{equation}
\end{proposition}
The case $n=m=0$ is the usual duality relation for multiple zeta values.
\subsection{Connected sum}
For $\bk, \bl\in\II_0$, the connected sum $Z_{n,m}^{\text{D}}(\bk;\bl)$ (\cite{SY1}) is defined by
\[
Z_{n,m}^{\text{D}}(\bk;\bl)\coloneqq\sum_{\substack{n=n_0<n_1<\cdots<n_a \\ m=m_0<m_1<\cdots<m_b}}\frac{1}{\bn^{\bk}}\cdot C^{\text{D}}(n_a,m_b)\cdot\frac{1}{\bm^{\bl}}.
\]
Here, the connector is
\begin{equation}\label{con}
C^{\text{D}}(n_a,m_b)\coloneqq\frac{n_a!\cdot m_b!}{(n_a+m_b)!}.
\end{equation}
By definition, there is a symmetry
\begin{equation}\label{eq:duality2}
Z_{n,m}^{\text{D}}(\bk;\bl)=Z_{m,n}^{\text{D}}(\bl;\bk).
\end{equation}
\subsection{Transport relations}
For $\bk,\bl\in\II_0$, we have
\begin{align}
Z_{n,m}^{\text{D}}(\bk_{\uparrow};\bl)&=Z_{n,m}^{\text{D}}(\bk;\bl_{\to})\quad (\bk\neq\varnothing),\label{eq:duality3}\\
Z_{n,m}^{\text{D}}(\bk_{\to};\bl)&=Z_{n,m}^{\text{D}}(\bk;\bl_{\uparrow})\quad (\bl\neq\varnothing).\label{eq:duality4}
\end{align}
\subsection{Algorithm}\label{D-alg}
For each $Z_{n,m}^{\text{D}}(\bk;\bl)$, we set rules as follows:

\

\begin{enumerate}[(i)]
\item\label{sec5:1} If $\bk\in\II'$, then use the transport relation \eqref{eq:duality3}.
\item\label{sec5:2} If $\bk\in\II\setminus\II'$, then use the transport relation \eqref{eq:duality4}.
\item\label{sec5:3} If $\bk=\varnothing$, then stop.
\end{enumerate}

\

Start from $Z_{n,m}^{\text{D}}(\bk;\varnothing)$ \ ($\bk\in\II'$) and transport indices from left to right according to the above rules until the algorithm stops.
Then, we have
\begin{equation}\label{eq:duality5}
Z_{n,m}^{\text{D}}(\bk;\varnothing)=Z_{n,m}^{\text{D}}(\varnothing;\bk^{\dagger}).
\end{equation}
\subsection{Boundary conditions}
For $\bk\in\II'$,
\[
Z_{n,m}^{\text{D}}(\bk;\varnothing)=\zeta^{}_{n,m}(\bk)
\]
holds by definition.
Therefore, by the symmetry \eqref{eq:duality2}, we see that the identity \eqref{eq:duality5} is nothing but the duality \eqref{eq:duality1}.
\subsection{Example}
By applying the algorithm in \ref{D-alg} to $Z_{n,m}^{\text{D}}((3,2);\varnothing)$, we have
\begin{align*}
Z_{n,m}^{\text{D}}((3,2);\varnothing)&=Z_{n,m}^{\text{D}}(\varnothing_{\to\uparrow\uparrow\to\uparrow};\varnothing)\\
&\stackrel{\eqref{sec5:1}}{=}Z_{n,m}^{\text{D}}(\varnothing_{\to\uparrow\uparrow\to};\varnothing_{\to}) \ =Z_{n,m}^{\text{D}}((3,1);(1))&\\
&\stackrel{\eqref{sec5:2}}{=}Z_{n,m}^{\text{D}}(\varnothing_{\to\uparrow\uparrow};\varnothing_{\to\uparrow}) \ \, =Z_{n,m}^{\text{D}}((3);(2))&\\
&\stackrel{\eqref{sec5:1}}{=}Z_{n,m}^{\text{D}}(\varnothing_{\to\uparrow};\varnothing_{\to\uparrow\to}) \ =Z_{n,m}^{\text{D}}((2);(2,1))&\\
&\stackrel{\eqref{sec5:1}}{=}Z_{n,m}^{\text{D}}(\varnothing_{\to};\varnothing_{\to\uparrow\to\to}) =Z_{n,m}^{\text{D}}((1);(2,1,1))&\\
&\stackrel{\eqref{sec5:2}}{=}Z_{n,m}^{\text{D}}(\varnothing;\varnothing_{\to\uparrow\to\to\uparrow}) \; =Z_{n,m}^{\text{D}}(\varnothing;(2,1,2))&
\end{align*}
and this proves
\[
\zeta^{}_{n,m}(3,2)=\zeta^{}_{m,n}(2,1,2).
\]
\section{Hoffman's identity}
\begin{proposition}[{\cite[Theorem~4.2]{H3}}]
For $\bk\in I_a$ $(a\geq 1)$,
\begin{equation}\label{eq:HD1}
H_N^{\star}(\bk)\coloneqq\sum_{1\leq n_1\leq\cdots\leq n_a\leq N}\frac{(-1)^{n_a-1}}{\bn^{\bk}}\binom{N}{n_a}=\zeta_N^{\star}(\bk^{\vee}).
\end{equation}
\end{proposition}
\subsection{Connected sum}
For $\bk\in\II$ and $\bl\in\II_0$, the connected sum $Z_N^{\text{HD}}(\bk;\bl)$ (\cite{SY2}) is defined by
\[
Z_N^{\text{HD}}(\bk;\bl)\coloneqq\sum_{1\leq n_1\leq\cdots\leq n_a\leq m_1\leq \cdots \leq m_b\leq m_{b+1}=N}\frac{1}{\bn^{(\bk_{\downarrow})}}\cdot C^{\text{HD}}(n_a,m_1)\cdot\frac{1}{\bm^{\bl}}.
\]
Here, the connector is
\[
C^{\text{HD}}(n_a,m_1)\coloneqq(-1)^{n_a-1}\binom{m_1}{n_a}
\]
and the connecting relation is $n_a\leq m_1$.
\subsection{Transport relations}
For $\bk\in\II$ and $\bl\in\II_0$, we have
\begin{align}
Z_N^{\text{HD}}(\bk_{\uparrow};\bl)&=Z_N^{\text{HD}}(\bk;{}_{\leftarrow}\bl),\label{eq:HD2}\\
Z_N^{\text{HD}}(\bk_{\to};\bl)&=Z_N^{\text{HD}}(\bk;{}_{\uparrow}\bl)\quad (\bl\neq\varnothing).\label{eq:HD3}
\end{align}
\subsection{Algorithm}\label{HD-alg}
For each $Z_N^{\text{HD}}(\bk;\bl)$, we set rules as follows:

\

\begin{enumerate}[(i)]
\item\label{sec6:1} If $\bk\in\II'$, then use the transport relation \eqref{eq:HD2}.
\item\label{sec6:2} If $\bk\in\II\setminus\II'$ and $\bk\neq (1)$, then use the transport relation \eqref{eq:HD3}.
\item\label{sec6:3} If $\bk=(1)$, then stop.
\end{enumerate}

\

Start from $Z_N^{\text{HD}}(\bk_{\uparrow};\varnothing)$ ($\bk\in\II$) and transport indices from left to right according to the above rules until the algorithm stops.
Then, we have
\begin{equation}\label{eq:HD4}
Z_N^{\text{HD}}(\bk_{\uparrow};\varnothing)=Z_N^{\text{HD}}((1);\bk^{\vee}).
\end{equation}
\subsection{Boundary conditions}
For $\bk\in\II$, we have
\[
Z_N^{\text{HD}}(\bk_{\uparrow};\varnothing)=H_N^{\star}(\bk)
\]
and
\[
Z_N^{\text{HD}}((1);\bk)=\zeta^{\star}_N(\bk).
\]
Therefore, we see that the identity \eqref{eq:HD4} is nothing but Hoffman's identity \eqref{eq:HD1}.
\subsection{Example}
By applying the algorithm in \ref{HD-alg} to $Z_{N}^{\text{HD}}((3,3);\varnothing)$, we have
\begin{align*}
Z_{N}^{\text{HD}}((3,2)_{\uparrow};\varnothing)&=Z_{N}^{\text{HD}}((1)_{\uparrow\uparrow\to\uparrow\uparrow};\varnothing)\\
&\stackrel{\eqref{sec6:1}}{=}Z_{N}^{\text{HD}}((1)_{\uparrow\uparrow\to\uparrow};{}_{\leftarrow}\varnothing) \ \, \, =Z_{N}^{\text{HD}}((3,2);(1))\\
&\stackrel{\eqref{sec6:1}}{=}Z_{N}^{\text{HD}}((1)_{\uparrow\uparrow\to};{}_{\leftarrow\leftarrow}\varnothing) \ =Z_{N}^{\text{HD}}((3,1);(1,1))\\
&\stackrel{\eqref{sec6:2}}{=}Z_{N}^{\text{HD}}((1)_{\uparrow\uparrow};{}_{\uparrow\leftarrow\leftarrow}\varnothing) \ \, =Z_{N}^{\text{HD}}((3);(2,1))\\
&\stackrel{\eqref{sec6:1}}{=}Z_{N}^{\text{HD}}((1)_{\uparrow};{}_{\leftarrow\uparrow\leftarrow\leftarrow}\varnothing) \ =Z_{N}^{\text{HD}}((2);(1,2,1))\\
&\stackrel{\eqref{sec6:1}}{=}Z_{N}^{\text{HD}}((1);{}_{\leftarrow\leftarrow\uparrow\leftarrow\leftarrow}\varnothing)=Z_{N}^{\text{HD}}((1);(1,1,2,1))
\end{align*}
and this proves
\[
H_N^{\star}(3,2)=\zeta^{\star}_N(1,1,2,1).
\]
\section{Cyclic sum formula}\label{sec:CS}
\begin{proposition}[{\cite[Theorem~2.3]{HO}}]
For a cyclic equivalent class $\alpha$ of an element of $\II'$, we have
\begin{equation}\label{eq:CS1}
\sum_{\bk\in\alpha}\zeta(\bk_{\uparrow})=\sum_{\bk\in\alpha}\sum_{j=0}^{k_a-2}\zeta({}_{\{\uparrow\}^j\leftarrow}\bk_{\{\downarrow\}^{j}}),
\end{equation}
where $\bk=(k_1,\dots,k_a)$.
\end{proposition}
\begin{proposition}[{\cite{KO}}]
For a cyclic equivalent class $\alpha$ of an element of $\II$, we have
\begin{equation}\label{eq:CS1'}
\sum_{\bk\in\alpha}\left(\zeta^{}_{p-1}(\bk_{\uparrow})+\zeta^{}_{p-1}({}_{\uparrow}\bk^{\sigma})+\zeta^{}_{p-1}({}_{\leftarrow}\bk^{\sigma})\right)\equiv \sum_{\bk\in\alpha}\sum_{j=0}^{k_a-2}\zeta^{}_{p-1}({}_{\{\uparrow\}^j\leftarrow}\bk_{\{\downarrow\}^{j}})\pmod{p},
\end{equation}
where $\bk=(k_1,\dots,k_a)$ and $\bk^{\sigma}=(k_a,k_1,\dots, k_{a-1})$.
\end{proposition}
\subsection{Connected sum}
For $\bk\in\II$, the connected sum $Z_N^{\text{O}}(\bk)$ is defined by
\[
Z_N^{\text{O}}(\bk)\coloneqq\sum_{0<n_1<\cdots<n_a\leq N}\frac{1}{\bn^{({}_{\downarrow}\bk)}}\cdot C^{\text{O}}(n_1,n_a).
\]
Here, the connector discovered by Ohno is
\[
C^{\text{O}}(n_1,n_a)\coloneqq\frac{1}{n_a-n_1}.
\]
\subsection{Transport relations}
For $\bk\in\II$, we have
\begin{align}
Z_N^{\text{O}}(\bk_{\uparrow})&=Z_N^{\text{O}}({}_{\uparrow}\bk)-\zeta^{}_N(\bk_{\uparrow}),\label{eq:CS2}\\
Z_{\infty}^{\text{O}}(\bk_{\to})&=Z_{\infty}^{\text{O}}({}_{\leftarrow}\bk)+\zeta(\bk_{\uparrow}),\label{eq:CS3}\\
Z_{p-1}^{\text{O}}(\bk_{\to})&\equiv Z_{p-1}^{\text{O}}({}_{\leftarrow}\bk)+\zeta^{}_{p-1}(\bk_{\uparrow})+\zeta^{}_{p-1}({}_{\uparrow}\bk)+\zeta^{}_{p-1}({}_{\leftarrow}\bk)\pmod{p}.\label{eq:CS4}
\end{align}
\subsection{Algorithm}\label{C-alg}
For each $Z_{\infty}^{\text{O}}(\bk)$ with $\bk\neq (1)$, we set rules as follows:

\

\begin{enumerate}[(i)]
\item\label{sec7:1} If $\bk\in\II'$, then use the case $N\to\infty$ of the transport relation \eqref{eq:CS2}.
\item\label{sec7:2} If $\bk\in\II\setminus\II'$, then use the transport relation \eqref{eq:CS3}.
\end{enumerate}

\

Start from $Z_{\infty}^{\text{O}}({}_{\leftarrow}\bk)$ ($\bk\in\II'$) and use transport relations according to the above rules until the first value $Z_{\infty}^{\text{O}}({}_{\leftarrow}\bk)$ appears again.
Then, we have
\begin{multline*}
Z_{\infty}^{\text{O}}({}_{\leftarrow}\bk)=Z_{\infty}^{\text{O}}({}_{\leftarrow}\bk)+(\text{L.H.S of \eqref{eq:CS1} for the class $[\bk]$})\\
-(\text{R.H.S of \eqref{eq:CS1} for the class $[\bk]$}).
\end{multline*}
\subsection{Boundary conditions}
We need the fact that the connected sum $Z_{\infty}^{\text{O}}(\bk)$ for an index $\bk$ one of whose component is greater than $1$ converges instead of boundary conditions (see \cite[Theorem~3.1]{HO}).

For $\bk\in\II$, we can also prove the cyclic sum formula for finite multiple zeta values \eqref{eq:CS1'} by using the transport relation \eqref{eq:CS4} instead of using the transport relation \eqref{eq:CS3} in the above algorithm.
\subsection{Example}
For $\bk=(k_1,\dots,k_a)\in I_a$ $(a\geq 1)$, we set
\[
S(\bk)\coloneqq\sum_{j=0}^{k_a-2}\zeta({}_{\{\uparrow\}^j\leftarrow}\bk_{\{\downarrow\}^{j}}).
\]
Here, $S(\bk)=0$ when $k_a=1$.
By applying the algorithm in \ref{C-alg} to $Z_{\infty}^{\text{O}}(1,2,1,3)$, we have
\begin{align*}
&Z_{\infty}^{\text{O}}({}_{\leftarrow}(2,1,3))\\
&\stackrel{\eqref{sec7:1}}{=}Z_{\infty}^{\text{O}}({}_{\uparrow\leftarrow}(2,1,3)_{\downarrow})-\zeta({}_{\leftarrow}(2,1,3))\\
&\stackrel{\eqref{sec7:1}}{=}Z_{\infty}^{\text{O}}({}_{\uparrow\uparrow\leftarrow}(2,1,3)_{\downarrow\downarrow})-\zeta^{}({}_{\uparrow\leftarrow}(2,1,3)_{\downarrow})-\zeta({}_{\leftarrow}(2,1,3))\\
&\stackrel{\eqref{sec7:2}}{=}Z_{\infty}^{\text{O}}({}_{\leftarrow}(3,2,1))+\zeta((3,2,1)_{\uparrow})-S(2,1,3)\\
&\stackrel{\eqref{sec7:2}}{=}Z_{\infty}^{\text{O}}({}_{\leftarrow}(1,3,2))+\zeta((1,3,2)_{\uparrow})+\zeta((3,2,1)_{\uparrow})-S(2,1,3)\\
&\stackrel{\eqref{sec7:1}}{=}Z_{\infty}^{\text{O}}({}_{\uparrow\leftarrow}(1,3,2)_{\downarrow})-\zeta({}_{\leftarrow}(1,3,2))+\zeta((1,3,2)_{\uparrow})+\zeta((3,2,1)_{\uparrow})-S(2,1,3)\\
&\stackrel{\eqref{sec7:2}}{=}Z_{\infty}^{\text{O}}({}_{\leftarrow}(2,1,3))+\zeta((2,1,3)_{\uparrow})-S(1,3,2)+\zeta((1,3,2)_{\uparrow})+\zeta((3,2,1)_{\uparrow})-S(2,1,3).
\end{align*}
This proves
\[
\zeta((2,1,3)_{\uparrow})+\zeta((1,3,2)_{\uparrow})+\zeta((3,2,1)_{\uparrow})=S(2,1,3)+S(1,3,2)+S(3,2,1),
\]
that is,
\[
\zeta(2,1,4)+\zeta(1,3,3)+\zeta(3,2,2)=\zeta(1,2,1,3)+\zeta(2,2,1,2)+\zeta(1,1,3,2).
\]
\section{Hoffman's relation}
\begin{proposition}[{\cite[Theorem~5.1]{H1}}]
For $\bk=(k_1,\dots,k_a)\in I'_a$ $(a\geq 1)$, 
\begin{equation}\label{eq:H1}
\sum_{i=0}^{a-1}\zeta(\bk_{(i)},{}_{\uparrow}\bk^{(i)})=\sum_{i=1}^{a}\sum_{j=1}^{k_i-1}\zeta((\bk_{(i)})_{\{\downarrow\}^j},{}_{\{\uparrow\}^j\leftarrow}\bk^{(i)}).
\end{equation}
\end{proposition}
\subsection{Connected sum}
For $\bk\in I_a$, $\bl\in I_b$ ($a,b\geq 1$), the connected sum $Z^{\text{H}}(\bk;\bl)$ is defined by
\[
Z^{\text{H}}(\bk;\bl)\coloneqq\sum_{0<n_1<\cdots<n_a<m_1<\cdots<m_b}\frac{1}{\bn^{(\bk_{\downarrow})}}\cdot C^{\text{H}}(n_a,m_1)\cdot\frac{1}{\bm^{\bl}}.
\]
Here, the connector is
\[
C^{\text{H}}(n_a,m_1)\coloneqq\frac{1}{m_1-n_a}
\]
and the connecting relation is $n_a<m_1$.
\subsection{Transport relations}
For $\bk\in\II$, $\bl\in\II'\cup\{(1)\}$, we have
\begin{align}
Z^{\text{H}}(\bk_{\uparrow};\bl)&=Z^{\text{H}}(\bk;{}_{\uparrow}\bl)+\zeta(\bk,{}_{\uparrow}\bl),\label{eq:H2}\\
Z^{\text{H}}(\bk_{\to};\bl)&=Z^{\text{H}}(\bk;{}_{\leftarrow}\bl)+\zeta(\bk,{}_{\leftarrow}\bl)\quad (\bl\neq (1)).\label{eq:H3}
\end{align}
\subsection{Algorithm}\label{H-alg}
For each $Z^{\text{H}}(\bk;\bl)$, we set rules as follows:

\

\begin{enumerate}[(i)]
\item\label{sec8:1} If $\bk\in\II'$, then use the transport relation \eqref{eq:H2}.
\item\label{sec8:2} If $\bk\in\II\setminus\II'$ and $\bk\neq(1)$, then use the transport relation \eqref{eq:H3}.
\item\label{sec8:3} If $\bk=(1)$, then stop.
\end{enumerate}

\

Start from $Z^{\text{H}}(\bk;(1))$ ($\bk\in\II'$) and transport indices from left to right according to the above rules until the algorithm stops.
Then, we have
\begin{equation}\label{eq:H4}
Z^{\text{H}}(\bk;(1))=Z^{\text{H}}((1);\bk)+\sum_{i=1}^{a}\sum_{j=1}^{k_i-1}\zeta((\bk_{(i)})_{\{\downarrow\}^j},{}_{\{\uparrow\}^j\leftarrow}\bk^{(i)})+\sum_{i=1}^{a-1}\zeta(\bk_{(i)},{}_{\leftarrow}\bk^{(i)}).
\end{equation}
\subsection{Boundary conditions}
For $\bk\in\II'$, we have
\[
Z^{\text{H}}(\bk;(1))=\sum_{i=0}^{a-1}\zeta(\bk_{(i)},{}_{\uparrow}\bk^{(i)})+\sum_{i=0}^{a-1}\zeta(\bk_{(i)},{}_{\leftarrow}\bk^{(i)})
\]
and
\[
Z^{\text{H}}((1);\bk)=\zeta({}_{\leftarrow}\bk).
\]
Therefore, we see that the identity \eqref{eq:H4} is nothing but Hoffman's relation \eqref{eq:H1}.
\subsection{Example}
For $\bk=(k_1,\dots,k_a)\in I'_a$ $(a\geq 1)$ and $1\leq i\leq a$, we set
\[
H_i(\bk)\coloneqq\sum_{j=1}^{k_i-1}\zeta((\bk_{(i)})_{\{\downarrow\}^j},{}_{\{\uparrow\}^j\leftarrow}\bk^{(i)}).
\]
Here, $H_i(\bk)=0$ when $k_i=1$. 
By applying the algorithm in \ref{H-alg} to $Z^{\text{H}}((2,1,3);(1))$, we have
\begin{align*}
&Z^{\text{H}}((2,1,3);(1))\\
&\stackrel{\eqref{sec8:1}}{=}Z^{\text{H}}((2,1,3)_{\downarrow};{}_{\uparrow\leftarrow}\varnothing)+\zeta((2,1,3)_{\downarrow},{}_{\uparrow\leftarrow}\varnothing)\\
&\stackrel{\eqref{sec8:1}}{=}Z^{\text{H}}((2,1,3)_{\downarrow\downarrow};{}_{\uparrow\uparrow\leftarrow}\varnothing)+\zeta((2,1,3)_{\downarrow\downarrow},{}_{\uparrow\uparrow\leftarrow}\varnothing)+\zeta((2,1,3)_{\downarrow},{}_{\uparrow\leftarrow}\varnothing)\\
&\stackrel{\eqref{sec8:2}}{=}Z^{\text{H}}((2,1);{}_{\leftarrow}(3))+\zeta((2,1),{}_{\leftarrow}(3))+H_3(2,1,3)\\
&\stackrel{\eqref{sec8:2}}{=}Z^{\text{H}}((2);{}_{\leftarrow}(1,3))+\zeta((2),{}_{\leftarrow}(1,3))+\zeta((2,1),{}_{\leftarrow}(3))+H_3(2,1,3)\\
&\stackrel{\eqref{sec8:1}}{=}Z^{\text{H}}((2)_{\downarrow};{}_{\uparrow\leftarrow}(1,3))+\zeta((2)_{\downarrow},{}_{\uparrow\leftarrow}(1,3))+\zeta((2),{}_{\leftarrow}(1,3))+\zeta((2,1),{}_{\leftarrow}(3))+H_3(2,1,3)\\
&=Z^{\text{H}}((1);(2,1,3))+H_1(2,1,3)+\zeta((2),{}_{\leftarrow}(1,3))+\zeta((2,1),{}_{\leftarrow}(3))+H_3(2,1,3).
\end{align*}
On the other hand, the boundary conditions are
\begin{align*}
Z^{\text{H}}((2,1,3);(1))&=\zeta({}_{\uparrow}(2,1,3))+\zeta((2),{}_{\uparrow}(1,3))+\zeta((2,1),{}_{\uparrow}(3))\\
&\qquad +\zeta({}_{\leftarrow}(2,1,3))+\zeta((2),{}_{\leftarrow}(1,3))+\zeta((2,1),{}_{\leftarrow}(3))
\end{align*}
and $Z^{\text{H}}((1);(2,1,3))=\zeta({}_{\leftarrow}(2,1,3))$.
Therefore, we have
\[
\zeta({}_{\uparrow}(2,1,3))+\zeta((2),{}_{\uparrow}(1,3))+\zeta((2,1),{}_{\uparrow}(3))=H_1(2,1,3)+H_2(2,1,3)+H_3(2,1,3),
\]
that is,
\[
\zeta(3,1,3)+\zeta(2,2,3)+\zeta(2,1,4)=\zeta(1,2,1,3)+\zeta(2,1,2,2)+\zeta(2,1,1,3).
\]
\section{Some remarks}
It is known that there exist various generalizations of each connected sum.
For example,
\[
\sum_{\substack{0=n_0<n_1<\cdots<n_a \\ 0=m_0<m_1<\cdots<m_b \\ n_a+m_b=r_1<r_2<\cdots<r_c\leq N}}\frac{x_1^{n_1}x_2^{n_2-n_1}\cdots x_a^{n_a-n_{a-1}}\cdot y_1^{m_1}y_2^{m_2-m_1}\cdots y_b^{m_b-m_{b-1}}\cdot z_2^{r_2-r_1}\cdots z_c^{r_c-r_{c-1}}}{\bn^{(\bk_{\downarrow})}\bm^{(\bl_{\downarrow})}\br^{({}_{\downarrow}\bh)}}
\]
and
\[
\sum_{\substack{1=r_0\leq r_1<\cdots<r_c \\ r_c=n_1<\cdots<n_a\leq N \\ r_c=m_1<\cdots<m_b\leq N}}\frac{x_2^{n_2}\cdots x_a^{n_a}\cdot y_2^{m_2}\cdots y_b^{m_b}\cdot z_1^{r_1}\cdots z_c^{r_c}}{\bn^{({}_{\downarrow}\bk)}\bm^{({}_{\downarrow}\bl)}\br^{\bh}}
\]
give the shuffle product formula and the harmonic product formula for multiple polylogarithms, respectively.
Here, $x_1, \dots, x_a, y_1,\dots, y_b, z_1,\dots, z_c$ are suitable complex variables.

We can also obtain a simple proof of the Ohno relation proved by Ohno in \cite{Oh} by generalizing the connector \eqref{con} as
\[
\frac{[n_a;x][m_b;x]}{[n_a+m_b;x]} \qquad \left(|x|<1, \ [n;x]\coloneqq \prod_{i=1}^n(i-x)\right).
\] 
Furthermore, we can get a simple proof of the $q$-Ohno relation proved by Bradley in \cite{B} by generalizing the above connector as
\[
\frac{[n_a;x]_q[m_b;x]_q}{[n_a+m_b;x]_q} \qquad \left(0<q<1, \ [n;x]_q\coloneqq \prod_{i=1}^n([i]_q-q^ix), \quad [i]_q=\frac{1-q^i}{1-q}\right).
\]
See \cite{SY1} for details.

We can also prove other series identities by using this dynamic proof method: the cyclic sum formula for (finite) multiple zeta-star values (same as \S\ref{sec:CS}), Leshchiner's identity which is a generalization of the Ap\'ery--Markov identity (Ono--Seki; unpublished), Zhao's binomial identity \cite[Theorem~1.4]{Z} which implies the two-one formlua (Yamamoto; unpublished), the duality for MZVs of level 2 (Ono--Seki, Yamamoto; both unpublished), and the double Ohno relation (Hirose--Sato--Seki; \cite{HSS}), and so on.

\

Let's find new connectors!
\section*{Acknowledgments}
This article is based on the talk by the author at the workshop ``Various Aspects of Multiple Zeta Values'' held at RIMS, Kyoto Univ.~(November 18--22,  2019). 
The author thanks the organizer Prof.~Hidekazu Furusho for the kind invitation.
This work was supported by JSPS KAKENHI Grant Number JP18J00151.
The author sincerely thanks Dr.~Minoru Hirose, Dr.~Masataka Ono, Dr.~Nobuo Sato, and Prof.~Shuji Yamamoto for their helpful comments. He also would like to thank Prof.~Masanobu Kaneko and Prof.~Yasuo Ohno for reading the manuscript carefully.

\end{document}